# EFFECTS OF STRATIGRAPHY ON RESPONSE OF ENERGY PILES

Arash Saeidi-Rashk-Olia, Department of Civil Engineering, Kansas State University, Manhattan, KS, Phone: 785-532-5862, e-mail: saeidi@ksu.edu
Dunja Perić, Department of Civil Engineering, Kansas State University, Manhattan, KS, Phone: 785-532-5862, e-mail: peric@ksu.edu

**ABSTRACT**

Energy piles are gaining increased popularity due to a growing demand for clean energy. To further advance the understanding of soil-structure interaction in energy piles, recently-derived analytical solutions have been implemented to investigate the impact of stratigraphy on the soil-structure interaction. This was accomplished by comparing the measured and predicted head displacements, axial strains and stresses in an energy pile embedded in the actual homogeneous and layered soil profiles, as well as into synthetic homogeneous soil profiles. The analytical solutions for homogeneous soil profile captured the smooth experimentally observed response versus depth very well. In the case of a layered soil profile, the corresponding analytical model was capable of capturing nuances in trends of axial stress and strain versus depth at the interface of different layers. The analytical predictions for the layered profile appear to be slightly less accurate than for the homogenous profile. The experimental data obtained from the layered profile appear to be a bit more scattered than those from the homogeneous profile. In the former case, the interplay of the individual soil layers with the pile occurs while maintaining the continuity of the pile stress and displacement at the interface of different layers. The response throughout each layer of the layered profile is quantitatively different than the corresponding homogeneous response. Nevertheless, qualitatively the response throughout each layer of the layered profile is similar to the response of the pile embedded in the corresponding homogeneous layer.

## INTRODUCTION

Use of geo-structures such as piles, walls, and tunnels couples their primary role as structural supports to the role of exchanging heat with shallow geothermal resources. This relatively new feature contributes to the sustainability of human constructions through the establishment of low-carbon buildings and infrastructure. Consequently, energy piles as sustainable foundations have a great potential for environmental and economic benefits, leading to their increasing popularity worldwide (Abdelaziz and Tolga 2016, Laloui and Rotta Loria 2019).
The free thermal deformation of energy piles is restrained by surrounding soil and rock layers, thus inducing thermal stress in the piles. Consequently, design of energy piles requires the knowledge of the thermally-induced stress and strain in the pile as well as of the interface shear and normal stresses (Elzeiny and Suleiman 2021). To accomplish this goal researchers have used physical and numerical models. Furthermore, analytical solutions have been found recently (Cossel 2019, Perić et al. 2020, Saeidi Rashk Olia and Perić 2021a).

One of the complexities in predicting the response of energy piles lies in the interaction of the pile shaft with various soil layers that may be present in the ground. The layers usually have different material properties, thus affecting the thermo-mechanical stresses that are induced in the piles. The current study focuses on implementing the recently-derived analytical solutions to advance the understanding of the impact of a site stratigraphy on the soil pile interaction in a single energy pile. To this end, the performance of analytical solutions for homogeneous and layered soil profiles is evaluated by comparing their predictions with the experimental results obtained from both, centrifuge tests and full-scale field tests. Furthermore, the effects of site stratigraphy and end restraints are also evaluated.

**ANALYTICAL SOLUTIONS**

Once an energy pile is subjected to heating, its tendency to expand induces an upward displacement in the upper part of the pile and downward displacement in its lower part. This leads to the emergence of a zero-displacement point known as the null point, which can be observed in Figures 1a and 1b. As shown, in both, homogeneous and layered soil configurations, the null point is located anywhere between the pile head and the pile tip. The heating of the pile induces compressive stress in the pile shaft, which is accompanied by downward pointing interface shear stress acting above the null point and upward-directed interface shear stress acting below the null point.

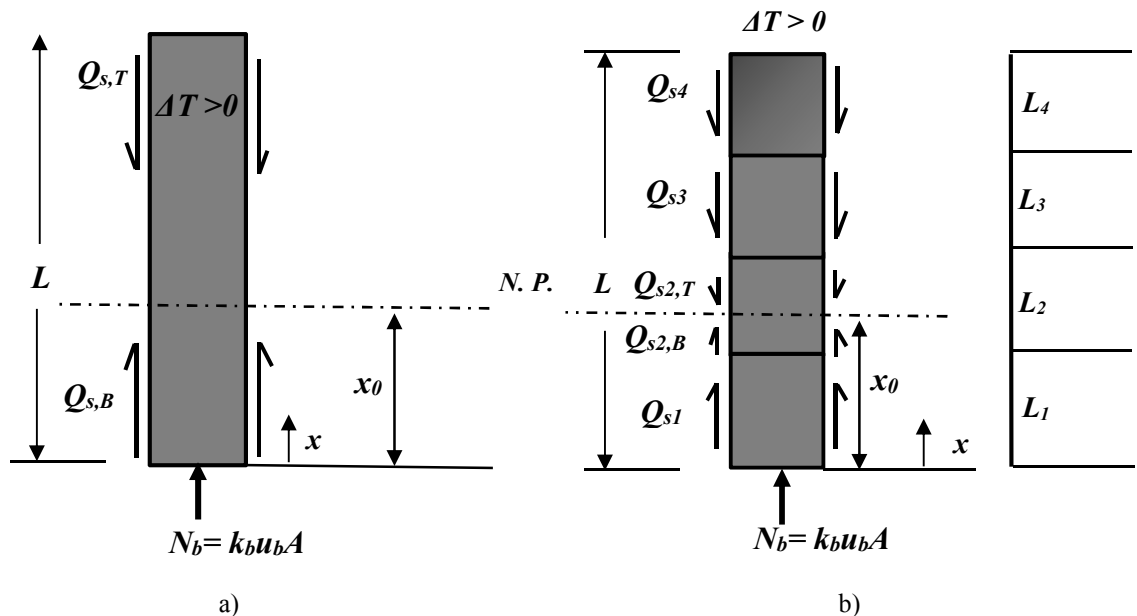

**Fig. 1. FBDs of the energy pile subjected to thermal load $\Delta T$ (a) in a homogeneous soil layer (b) in a four-layered soil profile**

Cossel (2019), Perić et al. (2020), and Saeidi Rashk Olia and Perić (2021a) derived analytical solutions for axial displacement, stress and strain in an energy pile subjected to thermo-mechanical loads. These solutions were derived based on the assumption of elastic soil pile interface behavior and thermo-elastic constitutive law for piles. These assumptions were previously confirmed by Knellwolf et al. (2011) and Perić et al. (2017). Cossel (2019) used a 1D formulation whereby the corresponding single coordinate was oriented in vertical direction. Both, the upward directed $x$-coordinate and displacement are positive as shown in Figure 1. Their solutions are based on the kinematic relationship of the pile given by

$$\varepsilon = \frac{du}{dx} \qquad [1]$$

where $u = u(x)$ denotes the axial displacement of the pile, and $\varepsilon = \varepsilon(x)$ is the axial strain of the pile. For a general case of thermo-elasticity, the corresponding constitutive law of the pile is given by

$$\sigma = E(\varepsilon - \alpha \Delta T) \qquad [2]$$

For mechanical loading $\Delta T = 0$, thus resulting in the Hooke's law. In these equations, $\sigma = \sigma(x)$ is the axial stress, $E$ is the Young's modulus of the pile, $\alpha$ is the coefficient of thermal expansion of the pile, and $\Delta T$ is a temperature change of the pile relative to the surrounding soil. The soil adjacent to the pile is modeled by a continuous linear elastic shear spring that is attached along the pile shaft. Thus, the shear stress $\tau = \tau(x)$, which acts along the pile shaft, was given by

$$|\tau| = k_s |u| \qquad [3]$$

where $k_s$ is the stiffness of the shear spring that is constant with depth within a given soil layer.

Enforcing the equilibrium of an infinitesimal pile segment, embedded in a homogeneous soil, in the vertical direction and using Eqs. 1, 2, and 3 results in

$$\frac{d^2 u}{dx^2} - \psi^2 u = 0 \qquad [4]$$

where the parameter $\psi$ is a constant given by

$$\psi^2 = \left(\frac{p}{A}\right)\left(\frac{k_s}{E}\right) = \xi_g \xi_s \qquad [5]$$

For a semi floating energy pile under thermo-mechanical loading, the boundary conditions at the pile tip and head are given by

$$\sigma(0) = k_b u(0) \qquad [6]$$

and

$$\sigma(L) = \frac{F}{A} \qquad [7]$$

where $k_b$ is the stiffness of the spring attached to the pile tip replicating the layers below the pile tip. By enforcing the boundary conditions from Eqs. 6 and 7 to the pile tip and head the solution for axial thermo-mechanical displacement, strain and stress in an energy pile is obtained. It is given by

$$u(x) = \overbrace{\frac{\alpha \Delta T \sinh[\psi(x-x_0)]}{\psi \cosh[\psi(L-x_0)]}}^{\text{Thermal}} + \overbrace{\frac{F[E\psi \cosh(\psi x) + k_b \sinh(\psi x)]}{AE\psi[E\psi \sinh(\psi L) + k_b \cosh(\psi L)]}}^{\text{Mechanical}} \qquad [8]$$

$$\varepsilon(x) = \frac{\alpha \Delta T \cosh[\psi(x-x_0)]}{\cosh[\psi(L-x_0)]} + \frac{F[E\psi \sinh(\psi x) + k_b \cosh(\psi x)]}{AE[E\psi \sinh(\psi L) + k_b \cosh(\psi L)]} \qquad [9]$$

$$\sigma(x) = E\alpha \Delta T \left[\frac{\cosh[\psi(x-x_0)]}{\cosh[\psi(L-x_0)]} - 1\right] + \frac{F[E\psi \sinh(\psi x) + k_b \cosh(\psi x)]}{A[E\psi \sinh(\psi L) + k_b \cosh(\psi L)]} \qquad [10]$$

whereby positive $x$ is directed upward, and $x_0$ is the location of the null point. Moreover, tensile strain and stress are positive. The location of thermal null point, $x_0$, is obtained by enforcing the boundary condition at the pile tip. It is given by

$$x_0 = \frac{1}{\psi} tanh^{-1}\left[\frac{cosh(\psi L)-1}{sinh(\psi L)+\frac{k_b}{E\psi}}\right] \quad [11]$$

Saeidi Rashk Olia and Perić (2021b and 2021c) implemented these solutions to investigate the effect of end restraints and mechanical axial load on thermo-mechanical response of energy piles during a complete heating-cooling cycle. Cossel (2019) and Perić et al. (2020) extended these solutions to a four-layered soil profile.

## PERFORMANCE OF ANALTYICAL MODEL FOR ENERGY PILE EMBEDDED IN A HOMOGENOUS SOIL PROFILE

Stewart and McCartney (2013) tested a concrete model of an end bearing single energy pile subjected to thermal load in a geotechnical centrifuge at a g-level of 24. The energy pile model had a diameter of 50.8 mm and a length of 533.4 mm, thus replicating a prototype energy pile having a length of 12.8 m and a diameter of 1.22 m. The energy pile had a modulus of elasticity of 7.17 GPa and it was embedded in a single layer of unsaturated Bonny silt with the coefficient of thermal expansion ($\alpha$) of $7.5 \times 10^{-6}$ 1/°C. While in real cases, the ideal zero displacements at the pile tip may never happen using this assumption could be reasonable for prediction of the load transfer behavior of piles resting on a very stiff bedrock. The experimental data suggest that this pile behaved as nearly as possible to an ideal end bearing pile exhibiting zero displacements at its tip. In order to incorporate the ideal end bearing behavior into the analytical solutions, the stiffness of the spring at the pile tip tends towards infinity, thus resulting in zero displacement at the pile tip. Consequently, this means that the null point is now located at the pile tip ($x_0 = 0$). After substituting $x_0 = 0$ and having $k_b$ tend to infinity, Eqs. 8 through 10 result in

$$u(x) = \frac{\alpha \Delta T sinh(\psi x)}{\psi\, cosh(\psi L)} + \frac{F sinh(\psi x)}{AE\psi cosh(\psi L)} \quad [12]$$

$$\varepsilon(x) = \frac{\alpha \Delta T cosh(\psi x)}{cosh(\psi L)} + \frac{F cosh(\psi x)}{AE cosh(\psi L)} \quad [13]$$

$$\sigma(x) = E\alpha\Delta T\left[\frac{cosh(\psi x)}{cosh(\psi L)} - 1\right] + \frac{F cosh(\psi x)}{A cosh(\psi L)} \quad [14]$$

Analytical predictions of the centrifuge tests are depicted in Figs. 2(a) through 2(c). The corresponding model parameters were obtained by validation of the analytical model that was performed by Saeidi Rashk Olia and Perić (2021). The resulting value of $k_s$ is equal to 55 MPa/m. Overall a very good fit between the analytical predictions and experimental data is observed over a range of temperature differences. These agreements confirm the validity of the end bearing analytical model derived with the assumption of ideal zero displacement at the pile tip that was incorporated into the analytical solution.

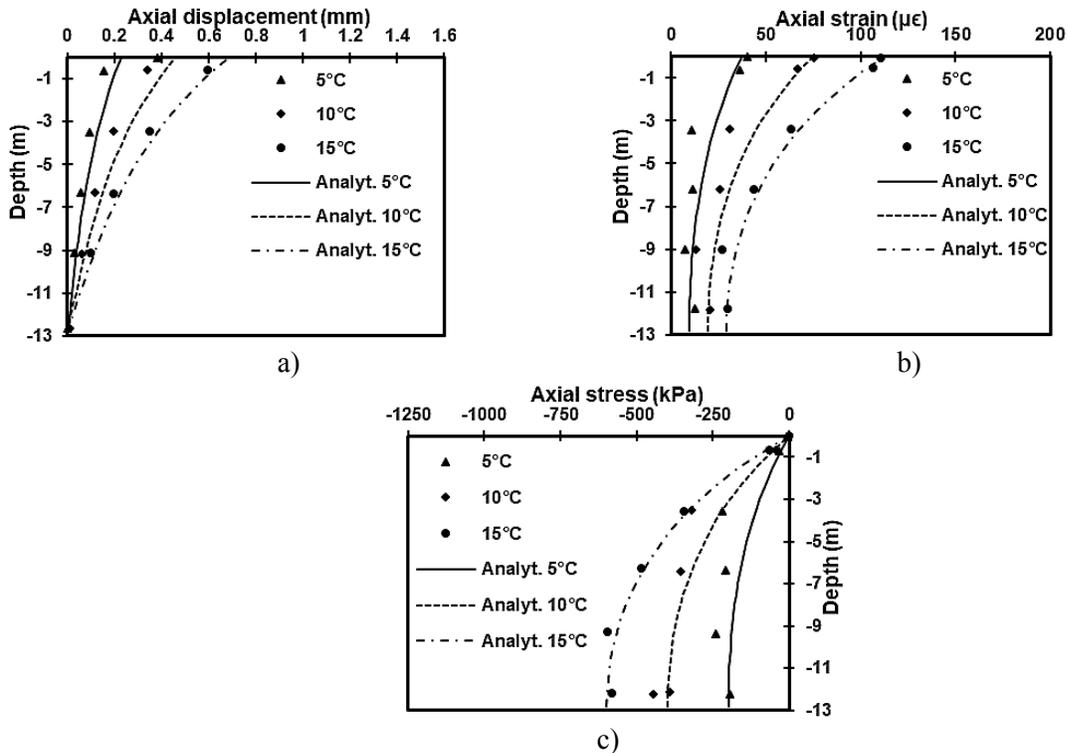

**Fig. 2. Comparison of analytical solutions and results of centrifuge tests for an end bearing energy pile: (a) thermal axial strain, (b) thermal axial stress, (c) thermal axial displacement**

The energy pile at the Swiss Federal Institute of Technology is embedded in the four-layered soil profile shown in Fig. 3. Perić et al. (2017) classified all four soil layers as low plasticity clays (CL) based on the laboratory test results reported by Laloui et al. (1999). Knellwolf et al. (2011) reported the stiffness of the elastic shear springs ($k_s$) to be equal to 16.7, 10.8, 18.2, and 121.4 MPa/m for soil layers A1, A2, B, and C, respectively based on the full-scale in situ pile tests. These values were used in the present study.

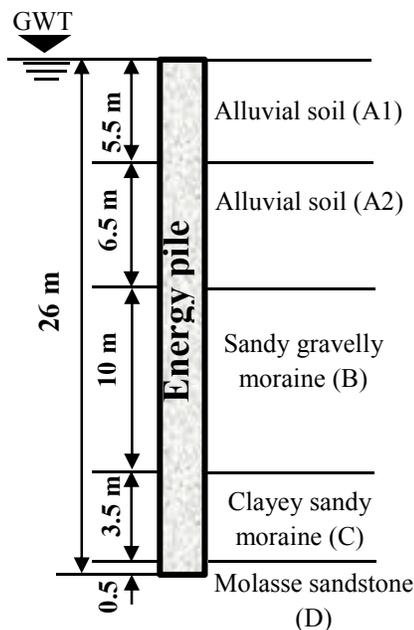

**Fig. 3. Layered soil profile at the Swiss Federal Institute of Technology, Lausanne**

Predictions of analytical solutions for the energy pile tested at the Swiss Federal Institute of Technology against experimental data provided by Laloui et al. (2006) were obtained herein based on the analytical solution for energy pile embedded into a four-layered soil (Cossel 2019). They are shown in Fig. 4 to Fig. 8. These figures include predictions for the pile embedded in the actual layered soil profile as well as for the pile embedded in hypothetical homogeneous profiles containing only soil A1, A2, B or C. Knellwolf et al. (2011) pointed out that while no precise information for the value $k_b$ for layer D (molasse) is available, this layer is the stiffest stratum of the geological profile. They proposed the stiffness of this layer to be one to two times the stiffness of the moraine (layer C), thus leading to $k_b$ values ranging between 667.7 MPa/m and 1335.5 MPa/m. The analyses conducted herein showed that the best fit with experimental data from thermal tests was obtained for $k_b$ = 6675 MPa/m, confirming the validity of the assumption that layer D is stiffer than layer C, however with a higher order of magnitude.

In Fig. 4, analytical predictions for the head displacement of the energy pile experiencing a heating-cooling cycle for both, layered and homogeneous profiles are shown along with experimental results. While the analytical prediction for the layered profile has the best fit with the experimental data, the predictions for homogeneous profiles comprised by soils A1 and B are also very close to the prediction for the layered profile. On the contrary, the prediction for the homogenous profiles consisting of soils A2 and C result in larger and significantly smaller displacements than those for the layered profile. This is expected because soils A2 and C have the lowest and highest $k_s$ values, respectively. Nevertheless, the value of $k_s$ for soil C is considerably larger than for any of the other three soils, thus resulting in a significantly smaller head displacement as compared to the three remaining soils and layered profile. It is noted that the values of $k_s$ for soils A1 and B are the closest to each other and that these soils cover 15.5m of the total 26 m of the pile length.

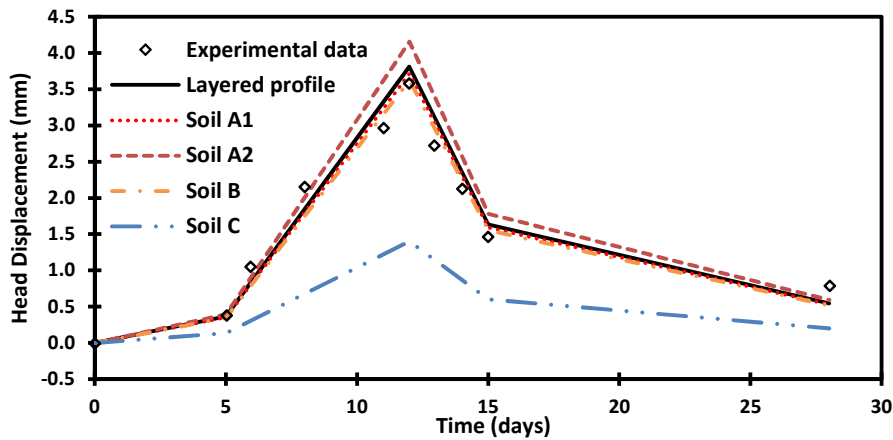

**Fig. 4. Analytical prediction of pile head displacement versus time for layered and homogeneous soil profiles**

Similar patterns of responses are observed in Fig. 5 and Fig. 6 that depict the axial thermal strain and axial stress for test T1, respectively. While the prediction for the layered profile provides the best fit with the experimental data overall, the predictions for the homogeneous profiles comprising soils A1 and B are not too far from the prediction for the layered profile. Furthermore, the prediction for the layered profile best fits the strain and stress measured at the pile tip. In addition, as expected, the homogeneous soil model is not able to capture the nuances in stress and strain response at the interface of different soil layers, especially at the boundary between the soils B and C, where the greatest change in $k_s$ occurs. Although the same value of $k_b$ was used for all homogeneous profiles as well as for the layered profile, the homogeneous profiles are not capable of matching the response at the pile tip due to their constant values of $k_s$ versus depth. Specifically, the predictions for the layered profile, homogeneous profiles containing soils A1 or B coincide throughout the depth of layer A1. Predictions for homogeneous profiles A1 and B start to deviate from the one for the layered profile upon entering into layer A2. Specifically, the strain predicted for homogeneous profiles containing soils A1 and B is smaller than the one predicted for the layered profile. On the contrary, the compressive stress predicted for homogeneous profiles containing soils A1 and B is larger than the one predicted for the layered profile. A slight increase in strain and decrease in stress is predicted by the model for the layered profile throughout layer B, as compared to layer A2. Finally, the model for layered profile predicts a further decrease in strain and increase in stress at the mid-depth of layer C as compared to the state at the interface of layers B and C.

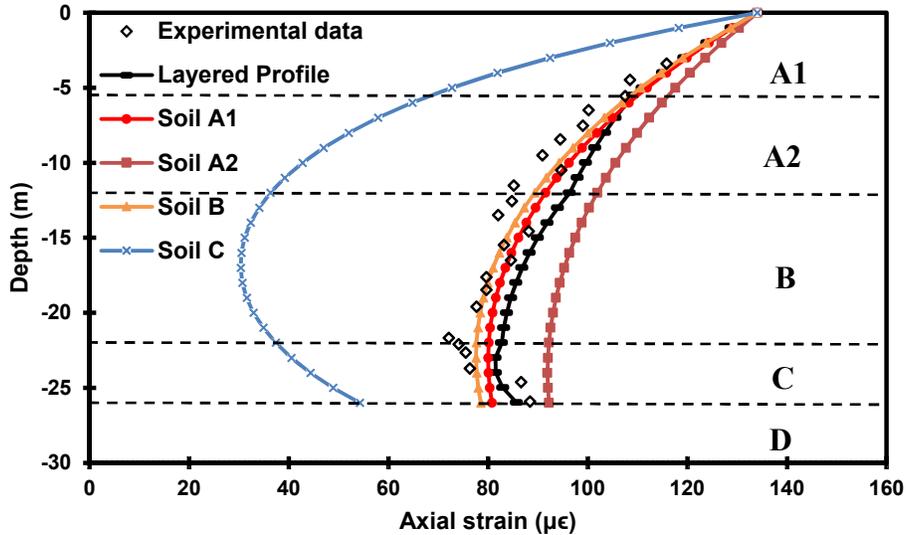

Fig. 5. Predictions of thermal axial strain versus depth for layered and homogeneous soil profiles (ΔT =13.4°C)

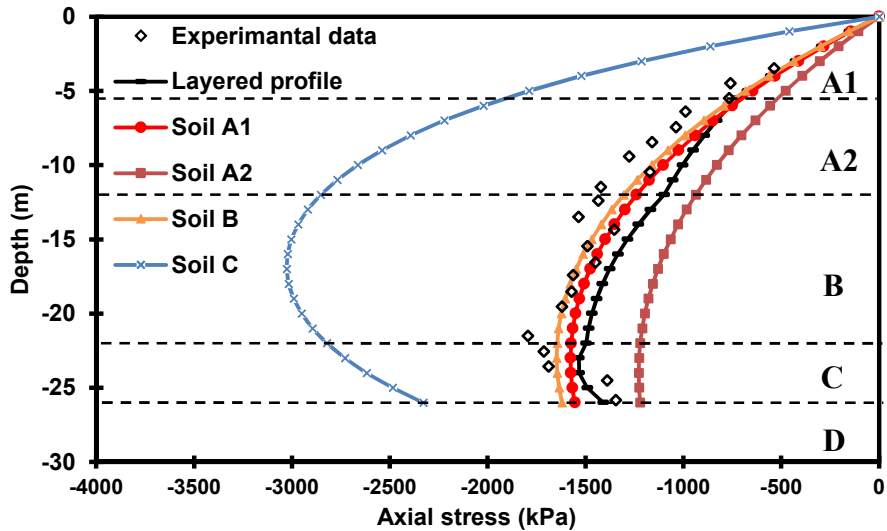

Fig. 6. Predictions of thermal axial stress versus depth for layered and homogeneous soil profiles (ΔT =13.4°C)

Predictions of axial strain and stress in the energy pile simultaneously subjected to thermal and mechanical loads (test T7) are presented in Fig.7 and Fig. 8. Similar trends in the prediction of the models for layered and homogeneous profiles are observed in the case of test T7 as in test T1 with the exception of the prediction of the model for a layered soil profile throughout the layer C. Specifically, the predicted strain is increasing, and stress is decreasing throughout layer C. As expected, the predictions of the model for the layered soil profile perform the best overall. They also qualitatively follow the experimental data, especially in the lower portion of the pile. Furthermore, it appears that the previously selected value of $k_b$ is almost perfect for test T1 but too large in the case of test T7. It is noted that both, the density and the scatter of the data obtained from the full-scale tests is larger than for the data obtained from the centrifuge tests.

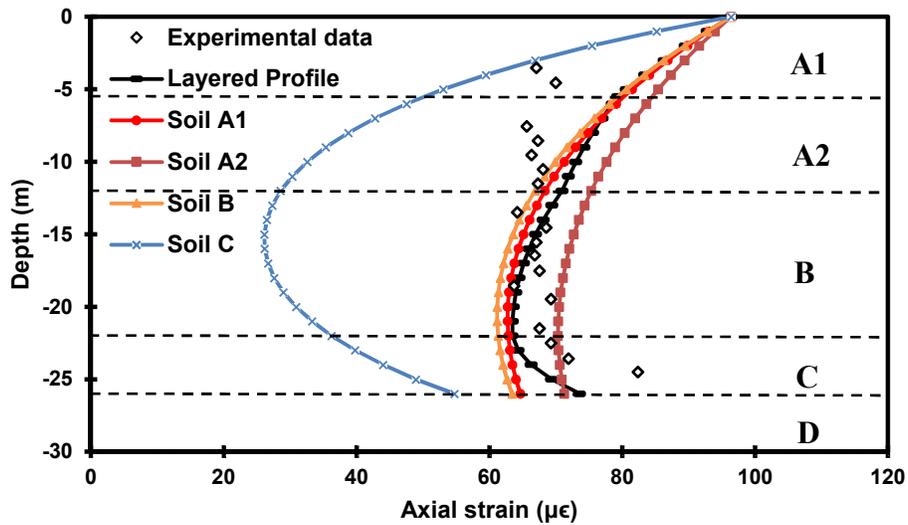

**Fig. 7. Prediction of thermo-mechanical axial strain versus depth for layered and homogeneous soil profiles (ΔT =14°C, F= -1000KN)**

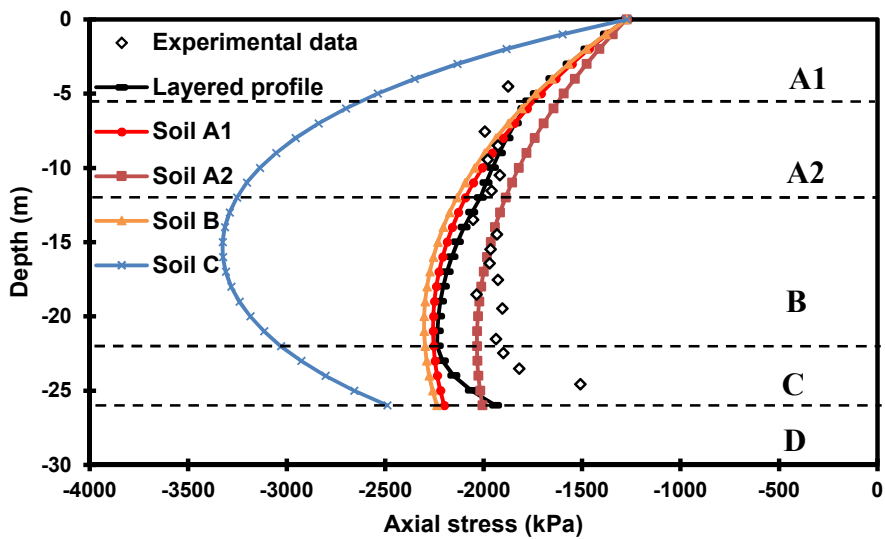

**Fig. 8. Prediction of thermo-mechanical axial stress versus depth for layered and homogeneous soil profiles (ΔT =14°C, F= -1000KN)**

# CONCLUSIONS

Recently-derived analytical solutions for thermal and mechanical responses of energy piles embedded in homogeneous and four-layered soil profiles were implemented to investigate the effects of stratigraphy on the soil-structure interaction in energy piles. To this end, analytical predictions for a total of six different soil profiles were compared to the experimental results of centrifuge tests and full scale field tests. The actual layered soil profile consisted of four different clay layers (A1, A2, B, C) located above the bedrock. For comparison, additional four synthetic homogeneous profiles were constructed, each containing a single soil (A1, A2, B, and C). Predictions of the analytical model for the actual homogeneous profile used in the centrifuge tests were slightly better than the predictions for the layered profile at the Lausanne site. Results indicate that the thermal and thermo-mechanical axial stress and strain distributions are dependent on the soil stratigraphy and predictions, since synthetic homogeneous soil layers are not completely capable of capturing the stress and strain variations along the entire length of the pile embedded in the layered soil profile. In particular, the stress and strain responses of the layered profile could exhibit an extremum value either at the boundary of the two layers or within one of the two layers. For the Lausanne site the former response was observed in the case of thermo-mechanical loading while the later was observed in the case of thermal loading. This occurred at the boundary of sandy gravelly moraine (layer B) and clayey sandy moraine (layer C). While in the case of the Lausanne site, some soil layers like A1 and A2 may be capable of predicting the overall trend of the pile response with good precision, where adjacent soil layers have largely different stiffness similar to soil layers B and C here, the homogeneity assumption can lead to underestimating the developed thermal in the pile.